\documentclass[10pt]{amsart}
\usepackage{amsmath,amstext,amssymb,amsopn,amsthm,mathrsfs}
\usepackage{graphicx}
\usepackage[english]{babel}
\usepackage{times}
\usepackage{fancybox}
\usepackage{epsfig}
\usepackage{graphics}
\usepackage{color}
\usepackage{amsmath}
\usepackage{amsfonts}
\usepackage[ansinew]{inputenc}
\usepackage{tikz}

\numberwithin{equation}{section}

\begin{document}

\title[Excursions on Cantor-like Sets] {Excursions on Cantor-like Sets}

\author{Robert DiMartino}
\address{Department of Mathematical and Actuarial Sciences, Roosevelt University  Chicago, Il, 60605, USA.}
\email{[Robert DiMartino]robert.dimartino@gmail.com}
\author{Wilfredo O. Urbina}
\email{[Wilfredo Urbina]wurbinaromero@roosevelt.edu}

\begin{abstract}
The ternary Cantor set $C$, constructed by George Cantor in 1883, is probably the best known example of a perfect nowhere-dense set in the real line, but as we will see later, it is not the only one. The present article we will explore the richness, the peculiarities and the generalities that has $C$ and explore some variants and generalizations of it. For a more systematic treatment the Cantor like sets we refer to our previous paper \cite{dimartinourb}.
\end{abstract}
\maketitle

\section{Introduction: Cantor Ternary Set.}
The Cantor  ternary set $C$ iis probably the best known example of a {\em perfect nowhere-dense} set in the real line. It was constructed by George Cantor in 1883, see \cite{cantor}.\\

$C$ is obtained from the closed interval $[0,1]$ by a sequence of deletions of open intervals known as ''middle thirds".
We begin with the interval $[0,1]$, let us call it $C_0$, and remove the middle third, leaving us with leaving us with the union of two closed intervals of length $1/3$
$$C_1= \left[0,\frac{1}{3}\right]  \cup \left[\frac{2}{3}, 1\right] .$$ Now we remove the middle third from each of these intervals, leaving us with the union of four closed intervals of length $1/9$
\begin{equation*}
C_2=\left[0,\frac{1}{9}\right] \cup \left[\frac{2}{9},\frac{1}{3}\right] \cup\left[\frac{2}{3},\frac{7}{9}\right] \cup\left[\frac{8}{9},1\right] .
\end{equation*}
Then we remove the middle third of each of these intervals leaving us with eight intervals of length $1/27$,
\begin{equation*}
C_{3}= [0,\frac{1}{27}]\cup[\frac{2}{27}, \frac{1}{9}]\cup[\frac{2}{9}, \frac{7}{27}]\cup [ \frac{8}{27},\frac{1}{3}] \cup [\frac{2}{3},\frac{19}{27}]\cup[\frac{20}{27},\frac{7}{9}]\cup[\frac{8}{9},\frac{25}{27}]\cup[\frac{26}{27},1].
\end{equation*}

We continue this process inductively, then for each $n=1,2, 3\cdots $ we get a set $C_n$  which is the union of $2^n$ closed intervals of length $1/3^n$. This iterative construction is illustrated in the following figure, for the first four steps:
\begin{center}
\includegraphics[width=3in]{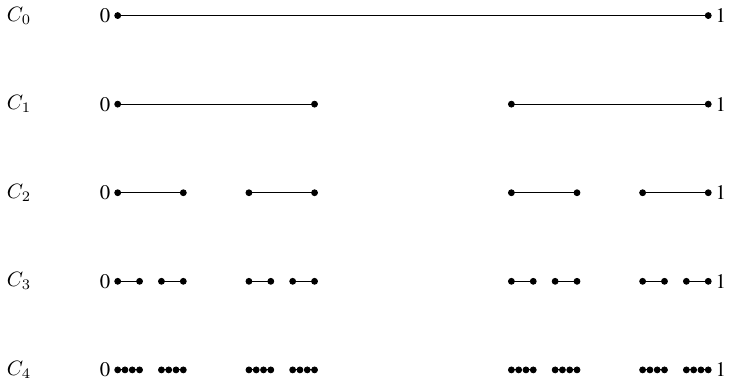}
\end{center}
Finally, we define the {\em Cantor ternary set} $C$ as the intersection
\begin{equation}
C= \bigcap_{n=0}^\infty C_n.
\end{equation}

Clearly $C \neq \emptyset $, since  trivially $0,1 \in C$. Moreover  $C$ is a a closed set, being the countable intersection of closed sets, and  trivially bounded, since it is a subset of $[0,1]$. Therefore, by the Heine-Borel theorem $C$ is a {\em compact set}. Moreover, observe by the construction that if $y$ is the end point of some closed subinterval of a given  $C_n$ then it is also the end point of some of the subintervals of $C_{n+1}$. Because at each stage, endpoints are never removed, it follows that $y \in C_n$ for all $n$. Thus $C$ contains all the end points of all the intervals that make up each of the sets $C_n$ (or alternatively,  the endpoints to the intervals removed)  all of which are rational ternary numbers in $[0,1]$, i.e. numbers of the form $k/3^n$. But $C$ contains much more than that; actually it is an uncountable set since it is a {\em perfect set} \footnote{A perfect set $P$ is a set that is closed and every point  $x \in P$ is a limit point i.e there is a sequence $\{x_n\} \subset P$, $x_n \neq x$ and $x_n \rightarrow x$.}. To prove that simply observe that every point of $C$ is approachable arbitrarily closely by the endpoints of the intervals removed (thus for any $x\in C$ and for each $n \in \mathbb{N}$ there is an endpoint, let us call it $y_n \in C_n$, such that $|x-y_n| < 1/3^n$).

There is an alternative characterization of $C$, the {\em ternary expansion characterization}.
Consider the ternary representation for  $x \in [0,1],$ \footnote{Observe, for the ternary rational  numbers $k/3^n$ there are two possible ternary expansions, since
$$ \frac{k}{3^n} = \frac{k-1}{3^n} + \frac{1}{3^n}  =  \frac{k-1}{3^n} + \sum_{k=n+1}^\infty \frac{2}{3^k}.$$ 
Similarly, for the dyadic rational numbers $k/2^n$ there are two possible dyadic expansions as
$$ \frac{k}{2^n} = \frac{k-1}{2^n} + \frac{1}{2^n}  =  \frac{k-1}{2^n} + \sum_{k=n+1}^\infty \frac{1}{2^k}.$$ 
Thus for the uniqueness of the dyadic and the ternary representations we will take  the infinite expansions representations for the dyadic and ternary rational numbers.}
 
\begin{equation}\label{ternaryexp}
 x = \sum_{k=1}^\infty \frac{\varepsilon_k(x)}{3^k}, \quad \varepsilon_k(x)=0, 1, 2 \quad \mbox{for all} \,k = 1, 2, \cdots 
\end{equation}

Observe that removing the elements where at least one of the \(\varepsilon_{k}\) is equal to one is the same as removing the middle third in the iterative construction, thus 
the Cantor ternary set is the set of numbers in $[0,1]$ that can be written in base 3 without using the digit 1, i.e.
\begin{equation}\label{ternaryChar}
 C=\left\{x\in [0,1] : x = \sum_{k=1}^\infty \frac{\varepsilon_k(x)}{3^k}, \quad \varepsilon_k(x)=0, 2 \quad \mbox{for all} \; k = 1, 2, \cdots
 \right\}.
\end{equation}

Using this characterization  of  $C$ we can get a direct proof that it is uncountable. Define 
the mapping $f: C \rightarrow [0,1]$ for $  x=\sum_{k=1}^\infty \frac{\varepsilon_k(x)}{3^k} \in C$, as 
\begin{equation}\label{cantorfun}
f(x) =    \sum_{k=1}^\infty \frac{\varepsilon_k(x)/2}{2^k}= \frac{1}{2}    \sum_{k=1}^\infty \frac{\varepsilon_k(x)}{2^k}.
\end{equation}

It is clear that $f$ is one-to-one correspondence from $C$ to $[0,1]$ (observe that as $\varepsilon_k=0, 2$ then $\varepsilon_k/2=0, 1$). 
As we have seen before, the un-contablity of $C$ can be also obtained from the fact that $C$ is perfect, see Abbott \cite{abot}, page 90. \\

$C$ is a {\em nowhere-dense} set, that is, there are no intervals included in $C$. One way to prove that is taking two arbitrary points in $C$ we can always find a number between them that requires the digit 1 in its ternary representation, and therefore there are no intervals included in $C$, thus $C$ is a nowhere dense set. Alternatively, we can prove this simply by contradiction. Assuming that there is a interval $I=[a,b] \subset C, \; a <b$. Then $I=[a,b] \subset C_n$ for all $n$ but as $|C_n| \rightarrow 0$ as $n \rightarrow \infty$ then $|I| = b-a =0.$\\

$C$ has measure zero, since its length can be obtained after subtracting from $1$ the sum of the length of all  open intervals removed in constructing it,
$$ m(C) = 1 - \sum_{n=1}^\infty \frac{2^{n-1}}{3^n} = 1 -\frac{1}{3}  \sum_{n=0}^\infty (\frac{2}{3})^{n} = 1 -  \frac{1/3}{1- 2/3}=1-1=0.$$\\

So far we have study the main properties of the Cantor set: $C$ is, non-empty (moreover uncountable) compact, perfect , nowhere dense set, with measure zero. Also we have seen that $C$ can then be obtained by any of the following constructions:
\begin{enumerate}
\item [i)] {\em Proportional (fractral) construction}: by removing a fixed proportion (one third) of each subinterval in each of the iterative steps.
\item[ii)] {\em Power construction}: by removing the length $1/3^n$ from the center of each subinterval in the $n^\text{th}-$step.
\item [iii)] {\em Ternary expansion characterization} : by removing  the digit one in each position of the ternary expansion.
\end{enumerate}

Of course these three constructions are equivalent. Nevertheless, the first construction, removing a fixed proportion of each subinterval in each of the iterative steps, give us another important property of the Cantor set, its self-similarity (fractal characteristic) across scales, this is illustrated in the following figure:
\begin{center}
\includegraphics[width=3in]{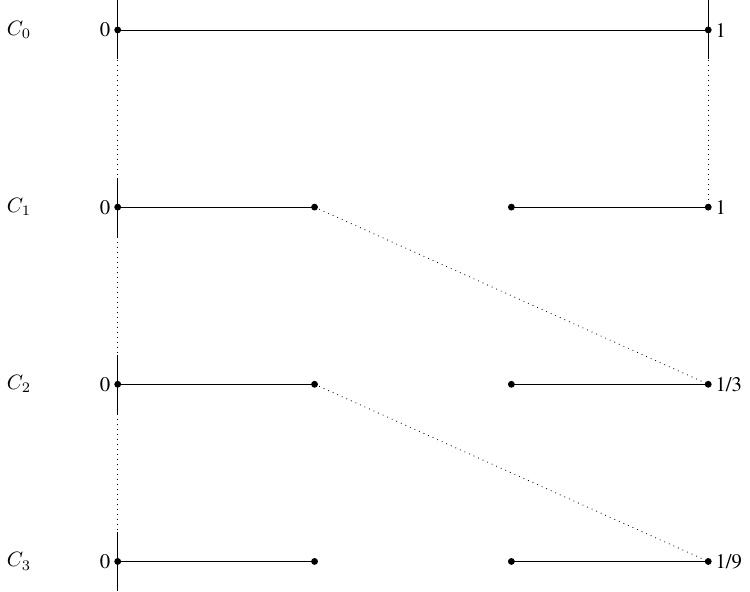}
\end{center}

In order to compute the Hausdorff dimension of $C$, S. Abbott in his book uses a nice little trick (avoiding the technicalities of the definition of Hausdorff dimension) by dilating $C_1$ by a the factor $3$ (for more details of this argument see \cite{abot} page 77) obtaining
$$ [0,1] \cup [2, 3].$$
Thus if we continue the iterations we will get two Cantor-like sets, then it can be concluded that
$$3^d= 2, $$
and therefore the Hausdorff dimension of $C$ is
 $$d= \frac{\ln 2}{\ln 3}=0.630930.$$
 
 The Cantor ternary set $C$ was not the first perfect nowhere-dense set  in the real line to be constructed. The first construction was done by the a British mathematician Henry J. S. Smith in 1875, but not many mathematicians were aware of Smith's construction. Vito Volterra, still a graduate student in Italy, also showed how to construct such a set  in 1881, but he published his result in an italian journal not widely read. In 1883 Cantor rediscover this construction himself, and due to Cantor's prestige, then the Cantor ternary set was  (and still is) the typical example of a perfect nowhere-dense set. Following D. Bresoud \cite{bres} we will refer as the Smith-Volterra-Cantor sets or $SVC(m)$ sets to the family of examples of perfect, nowhere-dense sets exemplified by the work of Smith, Voterra and Cantor, by removing in the $n$-iteration, an open interval of length $1/m^n$ from the center of the remaining closed intervals. Observe that $C = SVC(3)$.\\
 
Additionally, the Cantor set can also be obtained using an Iterated Function System (IFS), see \cite{hutch}. Consider the maps:
\begin{eqnarray}
\omega_0(x) = x/3, \quad \omega_1(x) = x/3 + 2/3,
\end{eqnarray}
then
$$ C = \omega_0([0,1]) \cup \omega_1([0,1]),$$
where this is a disjoint union. $C$ is the only set in the real line that satisfies this relation. This is a very interesting construction with important applications in the study of $C$, but we will not consider it in further details in this article.

\section{Variations of the Cantor ternary set:  Cantor-like sets.}

Is there something special with the number $3$? The answer is actually yes and no!. As we have already mentioned  $C$ can then be obtained by removing a fixed proportion (one third) of each subinterval in each of the iterative steps, by removing the length $1/3^n$ from each subinterval in the $n^\text{th}-$step, or removing  the digit one of the ternary expansion. What happen if we choose $2$ instead of $3$?
Using the proportional construction, removing the open ``middle half" from each component, we get at the end of the iterative process the set that we will denote as $C^{1/2}$.  Here are the first 3 iterative steps,

\begin{eqnarray*}
C^{1/2}_{1}&=& [0,\frac{1}{4}]\cup [\frac{3}{4},1],\\
C^{1/2}_{2}&=& [0,\frac{1}{16}]\cup[\frac{3}{16}, \frac{1}{4}]\cup [\frac{3}{4},\frac{13}{16}]\cup[\frac{15}{16},1],\\
\end{eqnarray*}
\begin{eqnarray*}
C^{1/2}_{3}&=& [0,\frac{1}{64}]\cup[\frac{3}{64}, \frac{1}{16}]\cup[\frac{3}{16}, \frac{13}{64}]\cup [ \frac{15}{64},\frac{1}{4}] \cup [\frac{3}{4},\frac{49}{64}]\cup[\frac{51}{64},\frac{13}{16}]\cup[\frac{15}{16},\frac{61}{64}]\cup[\frac{63}{64},1].
\end{eqnarray*}
Then, continuing this process inductively , then for each $n=1,2, \cdots $ we get a set $C^{1/2}_n$  which is the union of $2^n$ closed intervals of length $1/4^n$, and 
$$ C^{1/2} = \bigcap_{n=1}^\infty C^{1/2}_{n}.$$
$C^{1/2}$ share the same properties as $C$, i.e. it is a perfect, nowhere dense set in the real line. Also $C^{1/2}$ has measure zero, as  its length can be obtained after subtracting from $1$ the sum of the length of all  open intervals removed in constructing it,
\begin{eqnarray*}
1- \frac{1}{2} -2 \frac{1}{8}-4\frac{1}{32} +\cdots &=& 1- \frac{1}{2}[1+ \frac{1}{2}+\frac{1}{4} +\cdots ]\\
&=& 1- \frac{1}{2}  \sum_{n=0}^\infty (\frac{1}{2})^n =  1- \frac{1}{2}  \frac{1}{1-\frac{1}{2}} =  1- \frac{1}{2} \frac{1}{\frac{1}{2}} = 1- 1=0.
\end{eqnarray*}
In order to compute the Hausdorff dimension of this Cantor-like set, following the same argument as before, if we dilate $C_1^{1/2}$ by $4$ we obtain
$$ [0,1] \cup [3, 4].$$
Thus, if we continue the iterations we will get at the end two Cantor-like sets, therefore
$$4^d= 2,\; \mbox{so}\; d= \frac{\ln 2}{\ln 4}= \frac{\ln 2}{2\ln 2} = \frac{1}{2}= 0.5,$$
 so for this Cantor-like set the {\em dust} obtained  is {\em more sparse} than the one generated by $C$.
Also observe that, we can get a {\em expansion characterization} of $C^{1/2}$ but in base $n=4$, 
\begin{equation*}
 C^{1/2}=\left\{x\in [0,1] : x = \sum_{n=1}^\infty \frac{\varepsilon_n(x)}{4^n}, \quad \varepsilon_n(x)=0, 3 \quad \mbox{for all} \; n = 1, 2, \cdots
 \right\}.
\end{equation*}

On the other hand, what happen if we use the power construction? i.e.what happen if we remove the length $1/2^n$ from the center of each subinterval in the $n^\text{th}-$step?.  The first step is the same as before, obtaining the set
$$ [0,\frac{1}{4}]\cup [\frac{3}{4},1],$$
but next we need to remove two open intervals of length $1/2^2 = 1/4$, so we get only four points,
$$\{0, 1/4, 3/4,1\}$$
so the process stops there, at $n=2$, since there are no more intervals to remove form. Clearly this is not the same as removing  open ``middle half" from each component.\\

What happen if we consider now $n=4$? In this case both the proportional and the power construction can be iterated an infinity number of times, but, as we are going to see they turn out to be different sets.

\begin{enumerate}
\item [i)] {\em Proportional construction}: If we repeat the Cantor set's construction starting with the interval [0,1], removing the open ``middle fourth" from the center of the remaining closed intervals, we get at the end of the iterative process a set that we will denote as $ C^{1/4}$. Here are the first three iterations
\begin{eqnarray*}
C^{1/4}_1 &=& [0,\frac{3}{8}] \cup [\frac{5}{8}, 1]\\
C^{1/4}_2 &=& [0,\frac{9}{64}]\cup [\frac{15}{64}, \frac{3}{8}]\cup[\frac{5}{8},\frac{49}{64}]\cup[\frac{55}{64},1]\\
C^{1/4}_3 &=& [0,\frac{27}{512}]\cup [\frac{45}{512}, \frac{9}{64}]\cup[\frac{15}{64},\frac{147}{512}]\cup[\frac{165}{512},\frac{3}{8}]\cup[\frac{5}{8},\frac{347}{512}]\cup[\frac{365}{512},\frac{49}{64}]\cup[\frac{55}{64},\frac{467}{512}]\cup[\frac{485}{512},1].
\end{eqnarray*}
Then, continuing this process inductively, we get, for each $k=1,2, \cdots, $  a set $C^{1/4}_k$  which is the union of $2^k$ closed intervals of length $(3/8)^k$. 
$$ C^{1/4} = \bigcap_{n=1}^\infty C^{1/4}_{n}.$$
\begin{center}
\includegraphics[width=3in]{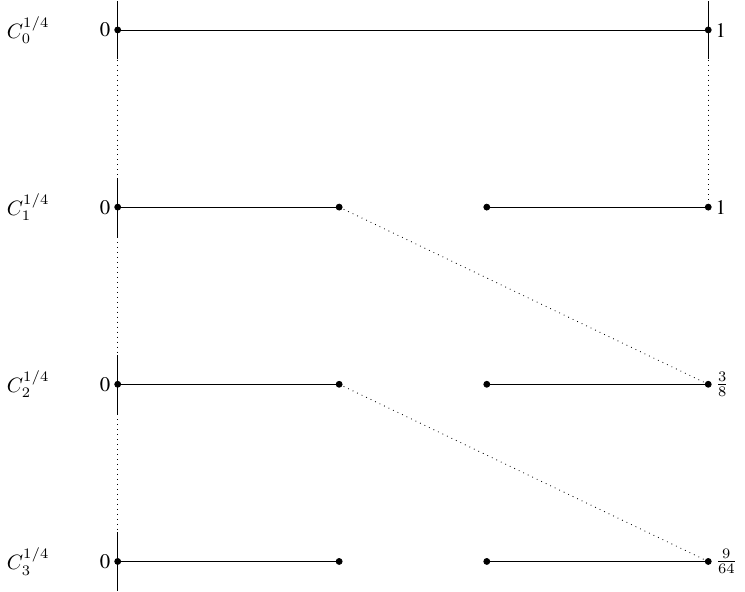}
\end{center}
$C^{1/4}$ share the same properties as $C$, i.e. it is a perfect, nowhere dense set in the real line.
Also $C^{1/4}$ has measure zero, as its length can be obtained after subtracting from $1$ the sum of the length of all  open intervals removed in constructing it,
\begin{eqnarray*}
1- \frac{1}{4} -2 \frac{3}{32}-4\frac{9}{256} +\cdots &=& 1- \frac{1}{4}[1+ \frac{3}{4}-\frac{9}{16} +\cdots ]\\
&=& 1- \frac{1}{4} \frac{1}{1-\frac{3}{4}} =  1- \frac{1}{4} \frac{1}{\frac{1}{4}} = 1-1 =0.
\end{eqnarray*}

In order to compute the Hausdorff dimension of this Cantor-like set, following the same argument as before, if we dilate $C^{1/4}_1$ by $8/3$ we obtain
$$ [0,1] \cup [\frac{5}{3}, \frac{8}{3}].$$
Thus if we continue the iterations we will get two Cantor-like sets, thus
$$(\frac{8}{3})^d= 2, \mbox{i.e} \quad d= \frac{\ln 2}{\ln 8/3} = \frac{\ln 2}{\ln 8- \ln3}= .706695.$$
Observe that if we do the same argument in the second or third step we get the same result,
$$(\frac{64}{9})^d= 4, \mbox{i.e} \quad d= \frac{\ln 4}{\ln 64/9} = \frac{\ln 2}{\ln 8- \ln3}= .706695,$$
$$(\frac{512}{27})^d= 8, \mbox{i.e} \quad d= \frac{\ln 8}{\ln 512/27} = \frac{\ln 2}{\ln 8- \ln3}= .706695.$$
Thus, the Cantor-like  set  $C^{1/4}$ is self-similar across scales (fractal)and it is less sparse than $C$.\\

Observe that $C^{1/4}$ can not be characterized using {\em expansion characterization}. It is easy to check that the ``natural base" $n=8$ does not work.\\

\item [ii)] {\em Power construction}: if we repeat the Cantor construction starting with the interval [0,1], removing, in the $n$-iteration, an open interval of length $1/4^n$ from the center of the remaining closed intervals, we get 
\begin{eqnarray*}
SVC(4)_1 &=& [0,\frac{3}{8}] \cup [\frac{5}{8}, 1]\\
SVC(4)_2 &=& [0,\frac{5}{32}]\cup [\frac{7}{32}, \frac{3}{8}]\cup[\frac{5}{8},\frac{25}{32}]\cup[\frac{27}{32},1]\\
SVC(4)_3 &=& [0,\frac{9}{128}]\cup [\frac{11}{128}, \frac{5}{32}]\cup[\frac{7}{32},\frac{37}{128}]\cup[\frac{39}{128},\frac{3}{8}]\cup[\frac{5}{8},\frac{89}{128}]\cup[\frac{91}{128},\frac{25}{32}]\cup[\frac{27}{32},\frac{59}{64}]\cup[\frac{119}{128},1].
\end{eqnarray*}
Then
$$ SVC(4) = \bigcap_{n=1}^\infty SVC(4) _{n}.$$
\begin{center}
\includegraphics[width=3in]{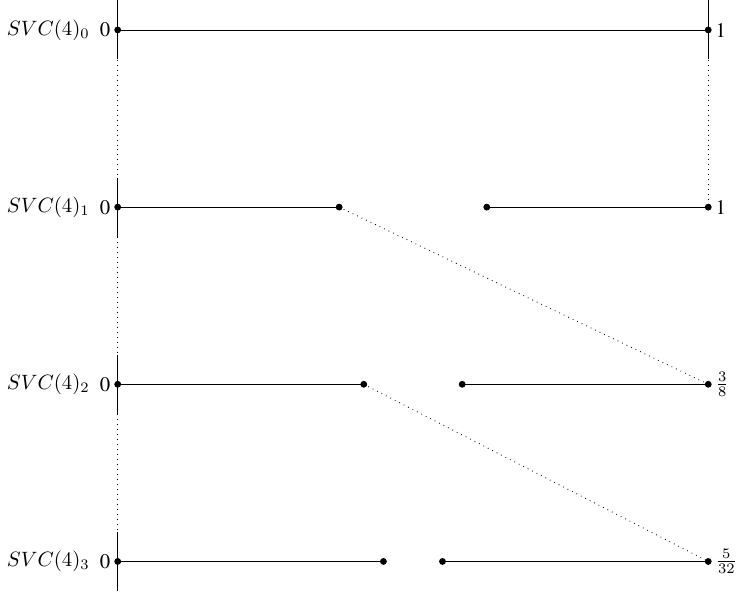}
\end{center}
 The total length of $SVC(4)$ can be obtained after subtracting from $1$ the sum of the length of all  open intervals removed in constructing it,
\begin{eqnarray*}
1- \frac{1}{4} -2 \frac{1}{16}-4\frac{1}{64} +\cdots &=& 1- \frac{1}{4}[1+ \frac{1}{2}-\frac{1}{4} +\cdots ]\\
&=& 1- \frac{1}{4} \frac{1}{1-\frac{1}{2}} =  1- \frac{1}{4} \frac{1}{\frac{1}{2}} = 1-\frac{1}{2} = \frac{1}{2}.
\end{eqnarray*}\\

Thus, the set $SVC(4)$ has positive measure equal to $1/2$. Cantor-like sets with positive measure are called {\em fat-Cantor sets}.\footnote{Observe that the only $SCV(n)$ set that has measure zero is $C$, since the total length of what is remove  to get it is
$$ 1 - \frac{1}{n} -\frac{2}{n^2} - \frac{4}{n^3}- \cdots = 1 - \frac{1}{n}[1 -\frac{2}{n} - \frac{2^2}{n^2}- \cdots ] =1 - \frac{1}{n} \frac{1}{1- 2/n} = \frac{n-3}{n-2}.$$}\\

In this case if we try to compute the Hausdorff dimension, since the sets $SVC(4)_n$ are not proportional,  a similar  argument as  before fails,  since if we dilate $SVC(4)_1$ by $8/3$ we obtain
$$ [0,1] \cup [\frac{5}{3}, \frac{8}{3}].$$
Thus if we continue the iterations we would get two Cantor-like sets, so
$$(\frac{8}{3})^d= 2, \mbox{i.e.} \quad d= \frac{\ln 2}{\ln 8/3} = \frac{\ln 2}{\ln 8- \ln3}= .706695.$$
But now if we do the same argument in the second or third step we get the different  results,
$$(\frac{32}{5})^d= 4, \mbox{i.e.} \quad d= \frac{\ln 4}{\ln 32/5} = \frac{2\ln 2}{\ln 32- \ln 5}= .746806,$$
$$(\frac{128}{9})^d= 8, \mbox{i.e.} \quad d= \frac{\ln 8}{\ln 128/9} = \frac{3\ln 2}{\ln 128- \ln 9}= .783274.$$
Moreover, since $$\lim_{k\rightarrow \infty} \frac{k \ln 2}{ \ln (2 \cdot 4^k) - \ln (1+ 2^k)} = 1,$$ 
we  can conclude that the Hausdorff dimension of $SVC(4)$  is actually one.\\
$SVC(4)$ is called the {\em Volterra set}. This is the set that was considered in 1881 by V.  Volterra to construct his famous counter-example of a function with bounded derivative that exists everywhere but the derivative is not Riemann integrable in any closed bounded interval i.e the Fundamental Theorem of Calculus fails!.\\

Observe that there is not a expansion characterization of the set $SVC(4)$.\\
\end{enumerate}

Now, if we repeat the Cantor set's {\em proportional construction} starting with the interval [0,1], removing the open ``middle three fourth" from the center of the remaining closed intervals, then we get in the first three iterations,
\begin{eqnarray*}
C^{3/4}_1 &=& [0,\frac{1}{8}] \cup [\frac{7}{8}, 1]\\
C^{3/4}_2 &=& [0,\frac{1}{64}]\cup [\frac{7}{64}, \frac{1}{8}]\cup[\frac{7}{8},\frac{57}{64}]\cup[\frac{63}{64},1]\\
C^{3/4}_3 &=& [0,\frac{1}{512}]\cup [\frac{7}{512}, \frac{1}{64}]\cup[\frac{7}{64},\frac{57}{512}]\cup[\frac{63}{512},\frac{1}{8}]\cup[\frac{7}{8},\frac{449}{512}]\cup[\frac{455}{512},\frac{57}{64}]\cup[\frac{63}{64},\frac{505}{512}]\cup[\frac{511}{512},1].
\end{eqnarray*}
Then, continuing this process inductively, for each $n=1,2, \cdots $ we get a set $C^{3/4}_kn$  which is the union of $2^n$ closed intervals of length $(3/8)^n$, and
$$ C^{3/4} = \bigcap_{n=1}^\infty C^{3/4}_{n}.$$
$C^{3/4}$ share the same properties as $C$, i.e. it is a perfect, nowhere dense set in the real line.
Also $C^{3/4}$ has measure zero, as 
\begin{eqnarray*}
1- \frac{3}{4} -2 \frac{3}{32}-4\frac{3}{256} +\cdots &=& 1- \frac{3}{4}[1+ \frac{1}{4}-\frac{1}{16} +\cdots ]\\
&=& 1- \frac{3}{4} \frac{1}{1-\frac{1}{4}} =  1- \frac{3}{4} \frac{1}{\frac{3}{4}} = 1-1 =0.
\end{eqnarray*}

In order to compute the Hausdorff dimension of this Cantor-like set, following the same argument as before, if we dilate $C^{3/4}_1$ by $8$ we obtain
$$ [0,1] \cup [7, 8].$$
Thus if we continue the iterations we will get two Cantor-like sets, thus
$$8^d= 2, \mbox{i.e} \quad d= \frac{\ln 2}{\ln 8} = \frac{\ln 2}{3\ln 2}= \frac{1}{3}= .666666,$$
 so this Cantor-like  set is less sparse than $C$. Also observe that in this case there is no {\em power construction.}\\
 
Finally, it is possible  an {\em expansion characterization} of that $C^{3/4}$ with base $8$,
\begin{equation}\label{Cant2}
 C^{3/4}=\left\{x\in [0,1] : x = \sum_{n=1}^\infty \frac{\varepsilon_n(x)}{8^n}, \quad \varepsilon_n(x)=0, 7 \quad \mbox{for all} \; n = 1, 2, \cdots
 \right\}.
\end{equation}
 \\

On the other hand, using the idea behind the ternary expansion characterization of the Cantor ternary set (\ref{ternaryChar}), we can get a general method to get Cantor-like sets. For instance, take the base $n=5$, and consider
$$C^5(0,2,4)= \left\{x\in [0,1] : x = \sum_{n=1}^\infty \frac{\varepsilon_n(x)}{5^n}, \quad \varepsilon_n(x)=0, 2, 4 \quad \mbox{for all} \; n = 1, 2, \cdots
 \right\}.$$

Since, for this set  we are removing all the numbers in $[0,1]$ such that in its $5$-adic expansion the $n$-th digit  $1$ or $3$, for each $n \geq 1$, then $C^5(0,2,4)$ is equivalent to the following {\em proportional construction}: start with the unit interval $[0,1]$, partition it into $5$  sub-intervals of equal length and remove the second and fourth open sub-intervals of the partition obtaining $C^5_{1}(0,2,4)$. 
$$C^5_{1}(0,2,4)= [0,\frac{1}{5}]\cup [\frac{2}{5},\frac{3}{5}]\cup [\frac{4}{5},1].$$
In the second iteration, partition each of the remain $3$ intervals in $C^5_{1}(0,2,4)$ into $5$ sub-intervals of equal length and  remove again second and fourth open sub-intervals, obtaining $C^5_{2}(0,2,4)$, 
\begin{eqnarray*}
C^5_{2}(0,2,4)&=& [0,\frac{1}{25}]\cup[\frac{2}{25}, \frac{3}{25}]\cup [\frac{4}{25},\frac{1}{5}]\cup[\frac{2}{5},\frac{11}{25}]\cup[\frac{12}{25},\frac{13}{25}]\cup[\frac{14}{25},\frac{3}{5}],\\
&& \quad \cup [\frac{4}{5},\frac{21}{25}]\cup[\frac{22}{25}, \frac{23}{25}]\cup [\frac{24}{25},1]\\
\end{eqnarray*}
One more iteration give us $C^5_{3}(0,2,4)$,
\begin{eqnarray*}
C^5_{3}(0,2,4)&=& [0,\frac{1}{125}]\cup[\frac{2}{125}, \frac{3}{125}]\cup [\frac{4}{125},\frac{1}{25}]\cup[\frac{2}{5},\frac{11}{125}]\cup[\frac{12}{125},\frac{13}{25}]\cup[\frac{14}{25},\frac{3}{25}],\\
&& \quad \cup [\frac{4}{25},\frac{21}{125}]\cup[\frac{22}{125}, \frac{23}{125}]\cup [\frac{24}{125},\frac{1}{5}]\ \cup [\frac{2}{5},\frac{51}{125}]\cup[\frac{52}{125}, \frac{53}{125}]\cup [\frac{54}{125},\frac{11}{25}]\\
&& \quad \cup [\frac{12}{25},\frac{61}{125}]\cup[\frac{62}{125}, \frac{63}{125}]\cup [\frac{64}{125},\frac{13}{125}]\ \cup [\frac{14}{25},\frac{71}{125}]\cup[\frac{72}{125}, \frac{73}{125}]\cup [\frac{74}{125},\frac{3}{5}]\\
&& \quad \cup [\frac{4}{5},\frac{101}{125}]\cup[\frac{102}{125}, \frac{103}{125}]\cup [\frac{104}{125},\frac{21}{25}]\ \cup [\frac{22}{25},\frac{111}{125}]\cup[\frac{112}{125}, \frac{113}{125}]\cup [\frac{114}{125},\frac{23}{52}]\\
&& \quad \cup [\frac{24}{25},\frac{121}{125}]\cup[\frac{122}{125}, \frac{123}{125}]\cup [\frac{124}{125},1].
\end{eqnarray*}
Continue partitioning each subinterval in $5$ sub-intervals and removing second and fourth open sub-intervals, we obtain for each $n\geq 1$,  $C^{5}_{n}$ which is the union of $3^{n-1}$ closed intervals of length $1/5^n$. Then
$$C^5(0,2,4)=\bigcap_{n=i}^{\infty} C^5_{n}(0,2,4).$$
By construction $ C^{5}(0,2,4)$ is also perfect nowhere set. Its length can be obtained after subtracting from $1$ the sum of the length of all  open intervals removed in constructing it,
\begin{eqnarray*}
1- 2\frac{1}{5} -6 \frac{1}{25}-18\frac{1}{125} +\cdots &=& 1- \frac{2}{5}[1+ \frac{3}{5}+\frac{9}{25} +\cdots ]\\
&=& 1- \frac{2}{5}  \sum_{k=0}^\infty (\frac{3}{5})^k = 1- \frac{2}{5}  \frac{1}{1-\frac{3}{5}} \\
&=&   1- \frac{2}{5} \frac{1}{\frac{2}{5}} = 1- 1=0.
\end{eqnarray*}
Thus $ C^{5}(0,2,4)$ must have measure $0$.\\

In order to compute the Hausdorff dimension of this Cantor-like set, following the same argument as before, if we dilate $C_1^{5}(0,2,4)$ by $5$ we obtain
$$ [0,1] \cup [2, 3] \cup [4, 5].$$
Thus if we continue the iterations we will get three Cantor-like sets, thus
$$5^d= 3, \; \mbox{i.e} \; d= \frac{\ln 3}{\ln 5}=0.6826,$$
so this Cantor-like  set is even less sparse than $C$.\\

Of course this construction can be generalized for any $m \in \mathbb{N}$ with $m>2$. Given $p= 1, \cdots m-1$ start again with the unit interval $[0,1]$.  Partition the interval into $m$  intervals of equal length and remove $p$ open intervals of the partition making sure to leave the first and last intervals\footnote{ if not, the set obtained would not be perfect since the construction would produce isolated points.} obtaining $C^m_{1}(p)$. In the second iteration, partition each of the remain $m-p$ intervals in $C^m_{1}(p)$ into $m$ sub-intervals of equal length and remove the corresponding $p$ open intervals to the ones removed in the first iteration, obtaining $C^m_{2}(p)$.  Continue partitioning each subinterval in $m$ sub-intervals and removing $p$ of them in the same pattern. Observe that, in this construction, in each iterative step $p$ subintervals are removed from each closed subinterval in  $C^m_{n}(p)$. Then
$$C^m(p)=\bigcap_{n=1}^{\infty} C^m_{n}(p).$$
\\
   
Clearly, the $m$-adic Cantor  set $C^m(p)$,  has the {\em expansion characterization} in base $m$
\begin{eqnarray*}
&& C^m(p) = \{ x \in [0,1] :  x=\sum_{n=1}^{\infty}\frac{\epsilon_{n}}{m^{n}}\,  \text{where }\epsilon_{n} \in \{ \eta_0, \eta_1, \cdots, \eta_p\} \subseteq \{0,1,\cdots,m-1\}, \\
&& \hspace{8cm} \,  \eta_0=0 \text{ and } \;\eta_p=m-1\} 
\end{eqnarray*}

Observe that in this case, using the dimension trick, we get
$$ m^d = p,$$
and therefore the Hausdorff dimension is
$$ d =\frac{\ln p}{\ln m}.$$\\

Several other constructions are possible. For instance,  take $ 0< \lambda \leq 1$ and repeat the construction of the Cantor set such that at the $n$-th step, instead of taking out an open interval of length $1/3^n$ from the center of each of the $2^{n-1}$ intervals of equal length, take out an open interval of length $\lambda/3^n$ to obtain $C_{\lambda,n}$, then, the first 3 iterative steps $C_{\lambda,1}$, $C_{\lambda,2}$, $C_{\lambda,3}$ are

\begin{eqnarray*}
C_{\lambda,1}&=& [0,\frac{3-\lambda}{6}]\cup [\frac{3+\lambda}{6},1],\\
C_{\lambda,2}&=& [0,\frac{9-5\lambda}{36}]\cup[\frac{9-\lambda}{36}, \frac{3-\lambda}{6}]\cup [\frac{3+\lambda}{6},\frac{27+\lambda}{36}]\cup[\frac{27+5\lambda}{36},1],\\
\end{eqnarray*}
\begin{eqnarray*}
C_{\lambda,3}&=& [0,\frac{27-19\lambda}{216}]\cup[\frac{27-11\lambda}{216}, \frac{9-5\lambda}{36}]\cup[\frac{9-\lambda}{36}, \frac{81-25\lambda}{216}]\cup [ \frac{81-17\lambda}{216},\frac{3-\lambda}{6}]\\
&& \cup [\frac{3+\lambda}{6},\frac{135+17\lambda}{216}]\cup[\frac{135+25\lambda}{216},\frac{27+\lambda}{36}]\cup[\frac{27+5\lambda}{36},\frac{189+11\lambda}{216}]\cup[\frac{189+19\lambda}{216},1].
\end{eqnarray*}
Then,
$$ C_{\lambda} = \bigcap_{n=1}^\infty C_{\lambda,n}$$
$ C_{\lambda}$ is, again, a perfect, nowhere dense set with positive measure, as
\begin{eqnarray*}
1- \frac{\lambda}{3} -2 \frac{\lambda}{9}-4\frac{\lambda}{27} +\cdots &=& 1- \frac{\lambda}{3}[1+ \frac{2}{3}-\frac{4}{9} +\cdots ]\\
&=& 1- \frac{\lambda}{3}  \sum_{k=0}^\infty (\frac{2}{3})^k =  1- \frac{\lambda}{3}  \frac{1}{1-\frac{2}{3}} =  1- \frac{\lambda}{3} \frac{1}{\frac{1}{3}} = 1- \lambda,
\end{eqnarray*}
Thus, for $ 0< \lambda< 1$,  $C_{\lambda}$ is a {\em fat Cantor-like} set. In order to compute the Hausdorff dimension, if we dilate $C_{\lambda,1}$ by $6/(3-\lambda)$ we obtain
$$ [0,1] \cup [\frac{3+\lambda}{3-\lambda}, \frac{6}{3-\lambda}].$$
Thus if we continue the iterations we will get two Cantor-like sets, thus
$$(\frac{6}{3-\lambda})^d= 2, \,\mbox{i.e.} \; d= \frac{\ln 2}{\ln [6/(3-\lambda)]} = \frac{\ln 2}{\ln 6- \ln(3-\lambda)}.$$
\vspace{0.25cm}

It is worth noting that, unlike the other previously explored constructions, in general Cantor-like sets  need not be symmetric. For instance,  $C^m(p)$ are not symmetric. As long as the first and last intervals of the {\em expansion characterization} partition remain, in order to avoid isolated points, any choice of removed open intervals can be made. Neither the measure nor the dimension of $C^m(p)$ is affected by the choice to remove intervals asymmetrically.\\

A more systematic study of the possible generalizations of the Cantor ternary set can be found in \cite{dimartinourb}.\\

Finally, a nice application of fat Cantor sets is the following. It is mention as one of the reasons the develop the theory Lebesgue integral the fact that the space of Riemann integrable functions is not a complete space. Unfortunately it is hard to find a concrete counterexample for that statement, i.e. a sequence of Riemann integrable functions such that its limit is not Riemann integrable. Using fat Cantor sets it is easy to construct such a counterexample. First of all remember that, according to a theorem due to H. Lebesgue, a function is Riemann integrable if and only if the set of discontinuities has measure zero, see for instance \cite{abot} Theorem 7.6.5 pag. 206. Let us take a fat Cantor set, for instance, let us take the Volterra set $SVC(4)$. We now know that $SVC(4)$ has measure $1/2>0$. Consider $K$ its complement in $[0,1]$ i.e. the union (disjunta) of the open intervals that have been removed in each iterative step,
$$K= [0,1] - SVC(4) = (\frac{3}{8}, \frac{5}{8}) \cup ( \frac{5}{32}, \frac{7}{32})\cup (\frac{25}{32},  \frac{27}{32})\cup ( \frac{9}{128}, \frac{11}{128}) \cup \cdots =\bigcup_{i=1}^\infty E_i.$$
Let  $f_n$ the indicator function of the set $\bigcup_{i=1}^n E_i,$ and $f$ be the indicator function of $K$. Then $f_n$ is Riemann integrable on $[0,1]$ since it is continuous for all but a finite number of points in [0,1]. Also trivially by construction $f_n (x)\rightarrow f(x),$ for any $x \in [0,1]$ and moreover, $f_n \rightarrow f$ in mean,
$$ \int_0^1 |f_n(x) - f(x) | dx = \int_{\bigcup_{i=n+1}^\infty E_i} dx = \sum_{i=n+1}^\infty |E_i| \rightarrow 0, $$
as $n \rightarrow \infty.$ But $f$ is not Riemann integrable on $[0,1]$ since its set of  discontinuities is precisely $SVC(4)$ that has positive measure.\\

\end{document}